\documentclass[10pt]{article}
\usepackage{amsmath, pb-diagram}
\usepackage{amssymb}
\usepackage{amscd}
\newenvironment{proof}{\par\noindent{\bf Proof \,}}{$\hfill \Box$\par\bigskip}
\date{\empty}
\newtheorem{thm}{Theorem}[section]
\newtheorem{lem}[thm]{Lemma}

\newtheorem{cor}[thm]{Corollary}

\newtheorem{ex}[thm]{Example}

\begin{document}

\title{Matrix nil-clean factorizations over abelian rings}

\author{Huanyin Chen\thanks{ Department of Mathematics, Hangzhou Normal University, Hangzhou,
310036, People's Republic of China, e-mail:
huanyinchen@aliyun.com}} \maketitle

\begin{abstract} A ring $R$ is nil-clean if
every element in $R$ is the sum of an idempotent and a nilpotent.
A ring $R$ is abelian if every idempotent is central. We prove
that if $R$ is abelian then $M_n(R)$ is nil-clean if and only if
$R/J(R)$ is Boolean and $M_n(J(R))$ is nil. This extend the main
results of Breaz et al. ~\cite{BGDT} and that of Ko\c{s}an et
al.~\cite{KLZ}.


\end{abstract}

\section{Introduction}
\vskip4mm Let $R$ be a ring with an identity. An element $a\in R$
is called nil-clean if there exists an idempotent $e\in R$ such
that $a-e\in R$ is a nilpotent. A ring $R$ is nil-clean provided
that every element in $R$ is nil-clean. In ~\cite[Question 3]{D},
Diesl asked: Let $R$ be a nil clean ring, and let $n$ be a
positive integer. Is $M_n(R)$ nil clean?  In ~\cite[Theorem
3]{BGDT}, Breaz et al. proved that their main theorem: for a field
$K$, $M_n(K)$ is nil-clean if and only if $K\cong {\Bbb Z}_2$.
They also asked if this result could be extended to division
rings. As a main result in ~\cite{KLZ}, Ko\c{s}an et al. gave a
positive answer to this problem. They showed that the preceding
equivalence holds for any division ring.

A ring $R$ is abelian if every idempotent in $R$ is central.
Clearly, every division ring is abelian. We extend, in this
article, the main results of Breaz et al.~\cite[Theorem 3]{BGDT}
and that of Ko\c{s}an et al.~\cite[Theorem 3]{KLZ}. We shall prove
that for an abelian ring $R$, $M_n(R)$ is nil-clean if and only if
$R/J(R)$ is Boolean and $M_n(J(R))$ is nil. As a corollary, we
also prove that the converse of a result of Ko\c{s}an et al.'s is
true.

Throughout, all rings are associative with an identity. $M_n(R)$
will denote the ring of all $n\times n$ full matrices over $R$
with an identity $I_n$. $GL_n(R)$ stands for the $n$-dimensional
general liner group of $R$.

\section{The main result}

We begin with several lemmas which will be needed in our proof of
the main result.

\begin{lem}~\cite[Proposition 3.15]{D} \label{21} Let $I$ be a nil ideal of a ring $R$. Then $R$ is nil-clean if and only if
$R/I$ is nil-clean.
\end{lem}

\begin{lem}~\cite[Theorem 3]{KLZ}\label{22} Let $R$ be a division ring. Then the following are equivalent:
\begin{enumerate}
\item [(1)]{\it $R\cong {\Bbb Z}_2$.}
\item [(2)]{\it $M_n(R)$ is nil-clean for all $n\in {\Bbb N}$.}
\item [(3)]{\it $M_n(R)$ is nil-clean for some $n\in {\Bbb N}$.}
\end{enumerate}
\end{lem}

Recall that a ring $R$ is an exchange ring if for every $a\in R$
there exists an idempotent $e\in aR$ such that $1-e\in (1-a)R$.
Clearly, every nil-clean ring is an exchange ring.

\begin{lem} \label{23} Let $R$ be an
abelian exchange ring, and let $x\in R$. Then $RxR=R$ if and only
if $x\in U(R)$.\end{lem} \begin{proof} If $x\in U(R)$, then
$RxR=R$. Conversely, assume that $RxR=R$. As in the proof of
~\cite[Proposition 17.1.9]{CH}, there exists an idempotent $e\in
R$ such that $e\in xR$ such that $ReR=R$. This implies that $e=1$.
Write $xy=1$. Then $yx=y(xy)x=(yx)^2$. Hence, $yx=y(yx)x$.
Therefore $1=x(yx)y=xy(yx)xy=yx$, and so $x\in U(R)$. This
completes the proof.\end{proof}

Set
$$J^*(R)=\bigcap \{ P~|~P~\mbox{is a maximal ideal of}~R\}.$$ We
will see that $J(R)\subseteq J^*(R)$. In general, they are not the
same. For instance, $J(R)=0$ and $J^*(R)=\{ x\in
R~|~dim_F(xV)<\infty\}$, where $R=End_F(V)$ and $V$ is an
infinite-dimensional vector space over a field $F$.

\begin{lem} \label{24} Let $R$ be an
abelian exchange ring. Then $J^*(R)=J(R)$.\end{lem} \begin{proof}
Let $M$ be a maximal ideal of $R$. If $J(R)\nsubseteq M$, then
$J(R)+M=R$. Write $x+y=1$ with $x\in J(R),y\in M$. Then $y=1-x\in
U(R)$, an absurd. Hence, $J(R)\subseteq M$. This implies that
$J(R)\subseteq J^*(R)$. Let $x\in J^*(R)$, and let $r\in R$. If
$R(1-xr)R\neq R$, then we can find a maximal ideal $M$ of $R$ such
that $R(1-xr)R\subseteq M$, and so $1-xr\in M$. It follows that
$1=xr+(1-xr)\in M$, which is imposable. Therefore $R(1-xr)R=R$. In
light of Lemma~\ref{23}, $1-xr\in U(R)$, and then $x\in J(R)$.
This completes the proof.\end{proof}

A ring $R$ is local if $R$ has only maximal right ideal. As is
well know, a ring $R$ is local if and only if for every $a\in R$,
either $a$ or $1-a$ is invertible if and only $R/J(R)$ is a
division ring.

\begin{lem}\label{25} Let $R$ be a ring with no non-trivial idempotents, and let $n\in {\Bbb N}$. Then the following are equivalent:
\begin{enumerate}
\item [(1)]{\it $M_n(R)$ is nil-clean.}
\item [(2)]{\it $R/J(R)\cong {\Bbb Z}_2$ and $M_n(J(R))$ is nil.}
\end{enumerate}
\end{lem}

\begin{proof} $(1)\Rightarrow (2)$ In view of ~\cite[Proposition 3.16]{D}, $J(M_n(R))$ is
nil, and then so is $M_n(J(R))$.

Let $a\in R$. By hypothesis, $M_n(R)$ is nil-clean. If $n=1$, then
$R$ is nil-clean. Then $a\in N(R)$ or $a-1\in N(R)$. This shows
that $a\in U(R)$ or $1-a\in U(R)$, and so $R$ is local. That is,
$R/J(R)$ is a division ring. As $R/J(R)$ is nil-clean, it follows
from Lemma~\ref{22} that $R/J(R)\cong {\Bbb Z}_2$. We now assume
that $n\geq 2$. Then there exists an idempotent $E\in M_n(R)$ and
a nilpotent $W\in GL_n(R)$ such that $I_n+\left(
\begin{array}{cccc}
a&&\\
&0&\\
&&\ddots&\\
&&&0
\end{array}
\right)=E+W$. Set $U=-I_n+W$. Then $U\in GL_n(R)$. Hence,
$$U^{-1}\left(
\begin{array}{cccc}
a&&\\
&0&\\
&&\ddots&\\
&&&0
\end{array}
\right)=U^{-1}E+I_n=\big(U^{-1}EU\big)U^{-1}+I_n.$$ Set
$F=U^{-1}EU$. Then $F=F^2\in M_n(R)$, and that
$$(I_n-F)U^{-1}\left(
\begin{array}{cccc}
a&&\\
&0&\\
&&\ddots&\\
&&&0
\end{array}
\right)=I_n-F.$$ Write $I_n-F=\left(
\begin{array}{cccc}
e&0&\\
*&0&\\
\vdots&&\ddots&\\
*&0&&0
\end{array}
\right).$ As $R$ possesses no non-trivial idempotents, $e=0$ or
$1$. If $e=0$, then $I_n-F=0$, and so $E=I_n$. This shows that
$\left(
\begin{array}{cccc}
a&&\\
&0&\\
&&\ddots&\\
&&&0
\end{array}
\right)=W$ is nilpotent; hence that $a\in R$ is nilpotent. Thus,
$1-a\in U(R)$.

If $e=1$, then $F=\left(
\begin{array}{cccc}
0&0&\\
*&1&\\
\vdots&&\ddots&\\
*&0&&1
\end{array}
\right).$ Write $U^{-1}=\left(
\begin{array}{cc}
\alpha&\beta\\
\gamma&\delta
\end{array}
\right),$ where $\alpha\in R, \beta\in M_{1\times (n-1)}(R),$
$\gamma\in M_{(n-1)\times 1}(R)$ and $\delta\in M_{(n-1)\times
(n-1)}(R)$. Then $$\left(
\begin{array}{cc}
\alpha&\beta\\
\gamma&\delta
\end{array}
\right)\left(
\begin{array}{cccc}
a&&\\
&0&\\
&&\ddots&\\
&&&0
\end{array}
\right)=\left(
\begin{array}{cc}
0&0\\
x&I_{n-1}
\end{array}
\right)\left(
\begin{array}{cc}
\alpha&\beta\\
\gamma&\delta
\end{array}
\right)+I_n,$$ where $x\in M_{(n-1)\times 1}(R)$. Thus, we get
$$\begin{array}{c}
\alpha a=1, \gamma a=x\alpha+\gamma, 0=x\beta+\delta+I_{n-1}.
\end{array}$$ One easily checks that
$$\left(
\begin{array}{cc}
1&\beta\\
0&I_{n-1}
\end{array}
\right)\left(
\begin{array}{cc}
1&0\\
x&I_{n-1}
\end{array}
\right)U^{-1}\left(
\begin{array}{cc}
1&0\\
\gamma a&I_{n-1}
\end{array}
\right)=\left(
\begin{array}{cc}
\alpha+\beta\gamma a&0\\
0&-I_{n-1}
\end{array}
\right).$$ This implies that $u:=\alpha+\beta\gamma a\in U(R)$.
Hence, $\alpha=u-\beta\gamma a$. It follows from $\alpha a=1$ that
$(u-\beta\gamma a)a=1$. As $R$ is connected, we see that
$a(u-\beta\gamma a)=1$, and so $a\in U(R)$. This shows that $a\in
U(R)$ or $1-a\in U(R)$. Therefore $R$ is local, and then $R/J(R)$
is a division ring. Since $M_n(R)$ is nil-clean, we see that so is
$M_n(R/J(R))$. In light of Lemma~\ref{22}, $R/J(R)\cong {\Bbb
Z}_2$, as desired.

$(2)\Rightarrow (1)$ In view of Lemma~\ref{21}, $M_n(R/J(R))$ is
 nil-clean. Since $M_n(R)/J(M_n(R))\cong M_n(R/J(R))$ and
$J\big(M_n(R)\big)=M_n(J(R))$ is nil, it follows from
Lemma~\ref{22} that $M_n(R)$ is nil-clean, as asserted.
\end{proof}

\begin{ex} \label{210} Let $K$ be a field, and let $R=K[x,y]/(x,y)^2$. Then $M_n(R)$ is nil-clean if and only if $K\cong {\Bbb Z}_2$.
Clearly, $J(R)=(x,y)/(x,y)^2$, and so $R/J(R)\cong K$. Thus, $R$ is a local ring with a nilpotent Jacobson radical. Hence, $R$ has no non-trivial idempotents.
Thus, we are done by Lemma~\ref{25}.\end{ex}

We come now to our main result.

\begin{thm}\label{26} Let $R$ be abelian, and let $n\in {\Bbb N}$. Then the following are equivalent:
\begin{enumerate}
\item [(1)]{\it $M_n(R)$ is nil-clean.}
\item [(2)]{\it $R/J(R)$ is Boolean and $M_n(J(R))$ is nil.}
\end{enumerate}
\end{thm}

\begin{proof} $(1)\Rightarrow (2)$ Clearly, $M_n(J(R))$ is nil. Let $M$ be a maximal ideal of
$R$, and let $\varphi_M: R\to R/M.$ Since $M_n(R)$ is nil-clean,
then so is $M_n(R/M)$. Hence, $R/M$ is an exchange ring with all
idempotents central. In view of ~\cite[Lemma 17.2.5]{CH}, $R/M$ is
local, and so $R/M$ is connected. In view of Lemma~\ref{25},
$R/M/J(R/M)\cong {\Bbb Z}_2$. Write $J(R/M)=K/M$. Then $K$ is a
maximal ideal of $R$, and that $M\subseteq K$. This implies that
$M=K$; hence, $R/M\cong {\Bbb Z}_2$. Construct a map $\varphi_M:
R/J^*(R)\to R/M, r+J^*(R)\mapsto r+M$. Then
$\bigcap\limits_{M}Ker\varphi_M=\bigcap\limits_{M}\{
r+J^*(R)~|~r\in M\}=0$. Therefore $R/J^*(R)$ is isomorphic to a
subdirect product of some ${\Bbb Z}_2$. Hence, $R/J^*(R)$ is
Boolean. In light of Lemma~\ref{24}, $R/J(R)$ is Boolean, as
desired.

$(2)\Rightarrow (1)$ Since $R/J(R)$ is Boolean, it follows by
~\cite[Corollary 6]{BGDT} that $M_n(R/J(R))$ is nil-clean. That
is, $M_n(R)/J(M_n(R))$ is nil-clean. But $J(M_n(R))=M_n(J(R))$ is
nil. Therefore we complete the proof, by
Lemma~\ref{21}.\end{proof}

We note that the "$(2)\Rightarrow (1)$" in Theorem ~\ref{26} always holds, but "abelian" condition is
necessary in "$(1)\Rightarrow (2)$". Let $R=M_n({\Bbb Z}_2) (n\geq 2)$. Then $R$ is
nil-clean. But $R/J(R)$ is not Boolean. Here, $R$ is not abelian.

\begin{cor}\label{27} Let $R$ be commutative, and let $n\in {\Bbb N}$. Then the following are equivalent:
\begin{enumerate}
\item [(1)]{\it $M_n(R)$ is nil-clean.}
\item [(2)]{\it $R/J(R)$ is Boolean and $J(R)$ is nil.}
\item [(3)]{\it For any $a\in R$, $a-a^2\in R$ is nilpotent.}
\end{enumerate}
\end{cor}

\begin{proof} $(1)\Rightarrow (3)$ Let $a\in R$. In view of Theorem~\ref{26},
$a-a^2\in J(R)$. Since $R$ is commutative, we see that $J(R)$ is
nil if and only if $J(M_n(R))$ is nil. Therefore $a-a^2\in R$ is
nilpotent.

$(3)\Rightarrow (2)$ Clearly, $R/J(R)$ is Boolean. For any $a\in
J(R)$, we have $(a-a^2)^n=0$ for some $n\geq 1$. Hence,
$a^n(1-a)^n=0$, and so $a^n=0$. This implies that $J(R)$ is nil.

$(2)\Rightarrow (1)$ As $R$ is commutative, we see that
$M_n(J(R))$ is nil. This completes the proof, by Theorem~\ref{26}.
\end{proof}

Furthermore, we observe that the converse of ~\cite[Corollary
7]{BGDT} is true as the following shows.

\begin{cor} \label{28} A commutative ring $R$ is nil-clean if and only if $M_n(R)$ is nil-clean.\end{cor} \begin{proof}
One direction is obvious by ~\cite[Corollary 7]{BGDT}. Suppose
that $M_n(R)$ is nil-clean. In view of Corollary~\ref{27} that
$R/J(R)\cong {\Bbb Z}_2$ is nil-clean, and that $J(R)$ is nil.
Therefore $R$ is nil-clean, by Lemma~\ref{21}.\end{proof}

\begin{ex} \label{211} Let $m,n\in {\Bbb N}$. Then $M_n\big({\Bbb Z}_m\big)$ is nil-clean if and only if $m=2^r$ for some $r\in {\Bbb N}$.
Write $m=p_1^{r_1}\cdots p_s^{r_s} (p_1,\cdots ,p_s~\mbox{are distinct primes}, r_1,\cdots ,r_s\in {\Bbb N}$).
Then $Z_m\cong Z_{p_1^{r_1}}\oplus \cdots \oplus Z_{p_m^{r_s}}$. In light of Corollary~\ref{27}, $M_n\big({\Bbb Z}_m\big)$ is nil-clean if and only if $s=1$ and
$Z_{p_1^{r_1}}$ is nil-clean. Therefore we are done
by Lemma~\ref{25}.\end{ex}

\end{document}